\documentclass[authoryear,review,3p]{elsarticle}

\usepackage{amsmath,amsthm,amssymb}
\usepackage[ruled,vlined]{algorithm2e}

\usepackage{booktabs}
\usepackage{makecell}

\usepackage{graphicx}
\usepackage{subfigure}
\usepackage{tikz}
\usepackage{pgfplots}
\pgfplotsset{compat=newest}
\usepgfplotslibrary{fillbetween}
\usetikzlibrary{3d,calc}

\usepackage[hidelinks]{hyperref}

\newcommand{\la}{\lambda}
\newcommand{\al}{\alpha}
\newcommand{\be}{\beta}

\newcommand{\R}{\mathbb{R}}

\newcommand{\norm}[1]{\left\Vert#1\right\Vert}

\newcommand{\ang}[1]{\left\langle{#1}\right\rangle}

\newtheorem{theorem}{Theorem}
\newtheorem{proposition}[theorem]{Proposition}
\newtheorem{lemma}[theorem]{Lemma}

\begin{document} 

\title{Cutting Plane Algorithms are Exact for Euclidean Max-Sum Problems}

\author[1,2]{Hoa T. Bui\corref{cor1}}
\ead{hoa.bui@curtin.edu.au}
\author[1,2]{Sandy Spiers}
\ead{sandy.spiers@postgrad.curtin.edu.au}
\author[1,2]{Ryan Loxton}
\ead{r.loxton@curtin.edu.au}
\address[1]{ARC Centre for Transforming Maintenance through Data Science, Curtin University, Perth, Australia}
\address[2]{Curtin Centre for Optimisation and Decision Science, Curtin University, Perth, Australia}
\cortext[cor1]{Corresponding author}

\begin{abstract}
    This paper studies binary quadratic programs in which the objective is defined by a Euclidean distance matrix, subject to a general polyhedral constraint set. 
    This class of nonconcave maximisation problems includes the capacitated, generalised and bi-level diversity problems as special cases.
    We introduce two exact cutting plane algorithms to solve this class of optimisation problems.
    The new algorithms remove the need for a concave reformulation, which is known to significantly slow down convergence.
    We establish exactness of the new algorithms by examining the concavity of the quadratic objective in a given direction, a concept we refer to as \emph{directional concavity}.
    Numerical results show that the algorithms outperform other exact methods for benchmark diversity problems (capacitated, generalised and bi-level), and can easily solve problems of up to three thousand variables.
\end{abstract}

\begin{keyword}
    Euclidean distance matrix \sep cutting plane methods \sep constrained diversity sum\sep exact algorithms \sep nonlinear binary optimisation 
\end{keyword}

\maketitle

\section{Introduction}
\label{sec:intro}

In this paper, we show that cutting plane algorithms are exact for the problem of maximising the sum of pairwise Euclidean distances between selected points, subject to general polyhedral constraints, hereafter referred to as the \emph{Euclidean max-sum problem}~\eqref{prob:emsp}.
The~\eqref{prob:emsp} is a generalisation of the Euclidean max-sum diversity problem \citep{spiers2023}, in which the cardinality constraint is replaced by a general polyhedral set.
More precisely, given a set of locations $v_1,\dots,v_n\in\R^s$ ($s\ge 1$), the~\eqref{prob:emsp} is defined as the following nonconcave binary maximisation problem,
\begin{align}
    \tag{EMSP}
    \label{prob:emsp}	
    \max\quad &f(x):=\tfrac{1}{2}\ang{Qx,x},\\
    \nonumber
    \text{s.t.}\quad & x\in P\cap \{0,1\}^n,
\end{align}
where $Q= [q_{ij}]_{i,j = 1,\ldots,n}$ is an $n\times n$ Euclidean distance matrix defined by $q_{ij} := \norm{v_i - v_j}$, and where $P\subset \mathbb{R}^n$ is a polyhedral set defined by
\begin{equation*}
    P = \left\{ x\in\R^n : Ax \le a \right\},
\end{equation*}
where $A \in \R^{m\times n}$ and $a\in\R^{m}$. 
Here, the definition of $x$ can be easily generalised to include both integer and continuous variables.
The matrix $Q$ is symmetric, hollow and has positive off-diagonal entries. 
By a result from~\cite{schoenberg1937certain}, given $v_1,\dots,v_n\in\R^s$, we can construct another set of $n$ points $u_1,\ldots,u_n\in\R^n$ such that $\norm{v_i - v_j}=\norm{u_i-u_j}^2$ for $i,j = 1,\ldots,n$. 
As such, the distance matrix $Q$ is also a \emph{squared} Euclidean distance matrix.
Furthermore, it is well-known that squared Euclidean distance matrices are \emph{conditionally negative definite}, i.e., $\ang{Qx,x} \le 0$ if $\sum_{i=1}^nx_i=0$, and have exactly one positive eigenvalue~\cite[Corollary 4.1.5, Theorem 4.1.7]{bapat_raghavan_1997}.
In this work, we exploit this property to prove that the cutting plane methodology, which is normally restricted to concave maximisation problems, converges to the optimal solution of~\eqref{prob:emsp}.

The Euclidean max-sum problem has various important practical applications.
In machine learning and statistical analysis, Euclidean distance is often used as a measure of dissimilarity between data points in clustering algorithms \citep{madhulatha2012overview,shirkhorshidi2015comparison}. 
By maximising the Euclidean distance between points, clusters can be formed based on their dissimilarity, allowing for effective grouping and classification of similar data.
An example of this is the well-known $k$-means clustering problem \citep{macqueen1967some_kmeans,lloyd1982least_kmeans}.
Furthermore, in various practical applications such as urban planning or network design, there is a need to strategically locate unwanted facilities such as waste disposal sites or polluting industries \citep{Kuby1987,ERKUT1989275}. 
Maximising the distance between these unwanted facilities and sensitive areas such as residential zones or environmental conservation areas helps minimise the negative impact on the surrounding communities or ecosystems.
Lastly, maximising Euclidean distances allows for the selection of points that capture diverse characteristics or represent different regions of interest, thereby enhancing the coverage and diversity of the chosen set.

This is seen in the Euclidean max-sum diversity problem \citep{spiers2023}, which is a special case of the~\eqref{prob:emsp} where the polyhedral set $P$ is defined by a single cardinality constraint.
For a recent review of this and other diversity models and their associated solution algorithms, we direct the reader to the comprehensive reviews in~\cite{Marti2022} and~\cite{parreno2021measuring}.
Among other applications, the maximum diversity problem has gained recent attention for its use in forming teams with diverse skill sets.

Recently, in~\cite{spiers2023}, we formulated a cutting plane algorithm for the Euclidean max-sum diversity problem by establishing the concavity of the objective function on the hyperplane $\sum_{i=1}^n x_i = p$, which ensures that tangent planes of feasible solutions serve as valid upper planes.
As such, our cutting plane algorithm is globally convergent for the Euclidean max-sum diversity problem.
The resultant exact algorithm is competitive with heuristic and meta-heuristic methods and can solve two coordinate instances of up to eighty thousand variables.
However, without the cardinality constraint, the objective function is not concave over the feasible set, and hence tangent planes do not always form valid cuts. 
The purpose of the current paper is to develop a new cutting plane methodology that still converges for this more general case, where concavity is not guaranteed.

To the best of our knowledge, outside of the Euclidean max-sum diversity problem, quadratic maximisation problems defined by Euclidean distance matrices have never been researched in isolation.
One reason for this is that these maximisation problems are, in general, nonlinear and nonconcave.
Mixed-integer nonlinear programming is one of the most challenging classes of optimisation problems.
Although there are several exact methods that provide general frameworks to tackle \emph{concave} maximisation problems, including outer approximation \citep{duran1986outer,leyffer1993deterministic,lubin2018polyhedral,kronqvist2020using}, branch and bound \citep{gupta1983nonlinear,vielma2008lifted,bonami2013branching}, and \emph{cutting plane methods} \citep{westerlund1995extended,kronqvist2016extended,lundell2022supporting}, advancements in exact algorithms for \emph{nonconcave} problems are still modest. The most common way to handle binary nonconcave maximisation is to reformulate the problem into an equivalent \emph{concave problem} by using a penalty approach, before applying exact methods to the new concave problem.
In particular, thanks to the property $x_i^2=x_i$ ($i=1,\ldots,n$) for $x\in\{0,1\}^n$, the nonconcave objective $f(x) = \tfrac{1}{2}\ang{Qx,x}$ can be replaced by a concave function ${f}_\rho(x) = \tfrac{1}{2}\ang{(Q-\rho I_n)x,x} + \tfrac{1}{2}\rho\sum_{i=1}^n x_i$, where $\rho$ is not smaller than the largest eigenvalue of the matrix $Q$. 
Although this technique is implemented in modern solvers such as \texttt{CPLEX} and \texttt{Gurobi} \citep{bliek1u2014solving,lima2017solution}, computational studies have shown that convergence is often slow \citep{lima2017solution,bliek1u2014solving,bonami2022classifier}. 
For the~\eqref{prob:emsp} where $Q$ is a Euclidean distance matrix, the Perron-Frobenius Theorem implies that the largest eigenvalue of $Q$ is bounded by the minimum and maximum row sums, and hence the concave reformulation requires choosing $\rho>0$.
We recently showed in Proposition 6,~\cite{spiers2023}, how choosing a large parameter $\rho$ in this concave reformulation step can weaken the cutting plane method. We further proved that if the polyhedral set $P$ is defined by a single cardinality constraint, then our cutting plane method converges to optimality without the need for concave reformulation.

This paper extends the results in~\cite{spiers2023} to general Euclidean distance maximisation by relaxing the requirement for a cardinality constraint.
This is achieved by exploiting the property that Euclidean distance matrices have exactly one positive eigenvalue.
To provide intuition to the reader on the key idea, consider a full eigenvalue decomposition of the objective function,
\begin{equation*}
    f(x) = \tfrac{1}{2}\ang{Qx,x} = \tfrac{1}{2} x^T \left( \sum_{i=1}^n \lambda_i v_i v_i^T \right) x = \tfrac{1}{2}\sum_{i=1}^n \lambda_i  x^T \left( v_iv_i^T \right) x, 
\end{equation*}
where $\{v_1,\ldots,v_n\}$ and $\la_1\ge \cdots \ge \la_n$ are eigenvectors and eigenvalues of $Q$.
This expresses the quadratic objective as a sum of functional components, which are either convex or concave depending on the sign of the respective eigenvalues.
However, as $Q$ is a Euclidean distance matrix, it is known to contain exactly one positive eigenvalue, and therefore $f(x)$ has only one convex component.
By restricting our search domain to exclude directions that traverse this convex component, our objective function can effectively be treated as a concave function (see Lemma~\ref{lemma:direction_concavity}).

\begin{figure}
    \centering
    \includegraphics[width=0.5\textwidth]{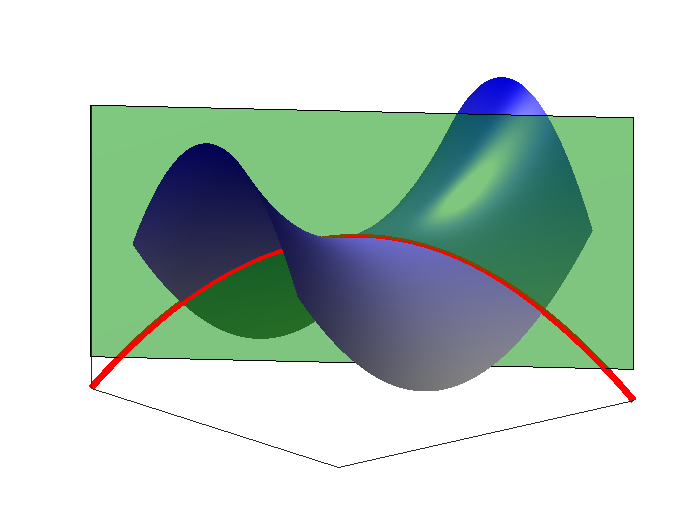}
    \caption{The intersection of a paraboloid and a hyperplane is either convex or concave.}
    \label{fig:paraboloid}
\end{figure}

To demonstrate this notion, consider the hyperbolic paraboloid defined by 
\begin{equation*}
    g(x,y) = xy = \tfrac{1}{4}(x+y)^2-\tfrac{1}{4}(x-y)^2.
\end{equation*}
Clearly, whenever $ax+by=0$ ($ab > 0$), the function reduces to $g(x,y)=x(- \tfrac{ax}{b}) = -\tfrac{a}{b}x^2$, which is concave.
Hence, while $g(x,y)$ is nonconcave for $x,y\in\R$, it is concave on the $ax+by=0$ plane.
The resultant concave parabola is shown in red in Figure~\ref{fig:paraboloid}.
This is essentially the technique used in~\cite{spiers2023}, where the Euclidean distance matrix is known to contain exactly one positive eigenvalue, and hence the objective has one convex functional component.
The cardinality constraint then ensures that the feasible set excludes this convex component, and the quadratic function can be treated as concave. 
For the general problem~\eqref{prob:emsp}, which may not include a cardinality constraint, the main idea of our approach is to only generate the tangent planes on concave directions. 
By doing so, the cutting planes are valid, and the algorithm always converges to an optimal solution.

The remainder of this paper is organised as follows.
In Section \ref{sec:method}, we formalise the concept of directional concavity and, based on this, formulate two key sufficient conditions for valid tangent planes, as detailed in Theorem~\ref{thm:valid_tangents}. 
These conditions then form the basis of two exact cutting plane algorithms, which vary in their approach to generating new cuts.
Finally, in Section \ref{sec:results} we conduct extensive numerical experiments to evaluate the effectiveness of the proposed solution approaches.

\section{Cutting plane algorithms}
\label{sec:method}

We denote the feasible set of~\eqref{prob:emsp} as $\mathcal{K}:=\left\{x\in \{0,1\}^n:\; x\in P\right\}\setminus \{0\}$. 
Note that we exclude $x = 0$ in $\mathcal{K}$ because $f(x) \ge 0 = f(0)$ for every $x\in \mathcal{K}$. 
Let $h: \R^n\times \R^n \to \R^n$ be the tangent plane of the function $f$, defined as:
\begin{equation*}
    h(x,y):= \ang{Qy, x-y} + \tfrac{1}{2}\ang{Qy,y}.
\end{equation*}
The tangent plane of a feasible solution $y\in\mathcal{K}$ is said to form a \emph{valid cut} if it provides an upper approximation for the optimal value of~\eqref{prob:emsp}, i.e, $f(x^*)\le h(x^*,y)$ where $x^*$ is an optimal solution of~\eqref{prob:emsp}. 

Since the function $f$ in~\eqref{prob:emsp} is not concave, not every feasible solution $y\in \mathcal{K}$ generates a valid cut. 
In this section, we establish sufficient conditions for when the tangent plane $h(x,y)$ is valid. 
The key to our approach is to study the concavity of the function $f$ when restricted to a given direction, exploiting the observation that the restriction of a quadratic function to a line is either concave or convex.

\subsection{Directional concavity}
\label{subsec:directional_concavity}

We explore the notion of \emph{directional concavity}; that is given a vector $u\in \R^n\setminus\{0\}$, we say that $u$ is a concave direction of $Q$ if $\ang{Qu,u} \le 0$.
Conversely, a vector $v\in \R^n\setminus\{0\}$ is a convex direction of $Q$ if $\ang{Qv,v} \ge 0$. 
Note that $x-y$ is a concave direction of $Q$ if and only if 
\begin{align*}
    h(x,y) - f(x) 
    &= \ang{Qy,x-y}+ \tfrac{1}{2}\ang{Qy,y} - \tfrac{1}{2}\ang{Qx,x} \\
    &= \ang{Qy,x-y}- \tfrac{1}{2}\ang{Q(x+y),x-y} \\
    &= -\tfrac{1}{2}\ang{Q(x-y),x-y}\ge 0.
\end{align*}
We first show that a vector $u = x-y$ is a concave direction of the matrix $Q$ if vector $u$ is orthogonal to $Qz$, where $z$ is a convex direction of $Q$.
\begin{lemma}
    \label{lemma:direction_concavity}
    Suppose $x,y \in \R^n$, and there is vector $z\in \R^n\setminus \{0\}$ such that
    \begin{enumerate}
        \item[a.] $\ang{Qz,z} \ge 0$, and
        \item[b.] $\ang{Qz,x-y} = 0$.
    \end{enumerate}
    Then, $h(x,y) \ge f(x)$.
\end{lemma}
\begin{proof}
	The inequality $h(x,y) \ge f(x)$ is equivalent to 
    \begin{equation*}
        \ang{Q(x-y),x-y} \le 0.
    \end{equation*}
	We suppose to the contrary that $\ang{Q(x-y),x-y} >0$. 
	Because $Q$ is a Euclidean Distance matrix, by~\cite[Corollary 4.1.5, Theorem 4.1.7]{bapat_raghavan_1997}, matrix $Q$ has exactly one positive eigenvalue. 
    Because $Q$ is a real symmetric matrix, it is orthogonally diagonalizable. 
    Let $\la_1> 0\ge \la_2\ge \ldots\ge \la_n$ be the eigenvalues of $Q$, and let $v_1,\ldots,v_n$ be the corresponding eigenvectors, which are normed, and orthogonal. 
    Then, we can express $x-y$ and $z$ on the basis $\{v_1,\ldots,v_n\}$ as follows:
    \begin{equation*}
        x-y = \sum_{i=1}^n \al_i v_i,\quad
        z =  \sum_{i=1}^n \be_i v_i,
    \end{equation*}
    for some $\al_i,\be_i\in \R$ ($i=1,\ldots,n$). 
    Then, 
    \begin{align}
        \label{proof:eq:z1} \ang{Qz,z} = \sum_{i=1}^n\la_i \be_i^2 \ge 0,\\
        \label{proof:eq:z2} \ang{Qz,x-y} = \sum_{i=1}^n\la_i\be_i\al_i = 0,\\
        \label{proof:eq:z3} \ang{Q(x-y),x-y} = \sum_{i=1}^n\la_i \al_i^2 >0.
    \end{align}
    Because $\la_i \le 0$ ($i=2,\ldots,n$), inequality~\eqref{proof:eq:z1} and $z\neq 0$ imply that $\beta_1 \neq 0$, and~\eqref{proof:eq:z3} implies that $\al_1 \neq 0$. 
    Therefore, we can multiply both sides of~\eqref{proof:eq:z1} by $\al_1^2 >0$, ~\eqref{proof:eq:z2} by $-2\al_1\be_1\neq 0$, and ~\eqref{proof:eq:z3} by $\be_1^2 >0$, and sum up to obtain
    \begin{align*}
        0 &< \left(\la_1 \be_1^2\al_1^2+\al_1^2\sum_{i=2}^n\la_i \be_i^2\right) - 2\left(\la_1 \be_1^2\al_1^2+\al_1\beta_1\sum_{i=2}^n\la_i \be_i\al_i\right)+ \left(\la_1 \be_1^2\al_1^2+\be_1^2\sum_{i=2}^n\la_i \al_i^2\right)\\
        &= \sum_{i=2}^n\la_i (\al_1^2\beta_i^2 - 2\al_1\be_1\al_i\be_i + \be_1^2\al_i^2)= \sum_{i=2}^n\la_i (\al_1\be_i - \al_i\be_1)^2.
    \end{align*}
    The inequality above only holds when there is at least one positive eigenvalue among $\la_2,\ldots,\la_n$, which is a contradiction. 
    Hence, it must hold that $\ang{Q(x-y), x-y} \le 0$.
\end{proof}

Recall that the Euclidean distance matrix $Q$ is \emph{conditionally negative definite}. 
The next result exploits this fact to replace condition (b) in Lemma~\ref{lemma:direction_concavity} with two new conditions.

\begin{lemma}
    \label{lemma:setup_valid_tangents}
    Suppose $x,y\in\R^n$, and there is $z\in\R^n\setminus\{0\}$ such that
    \begin{enumerate}
        \item[a.] $\ang{Qz,z}\ge0$,
        \item[b.] $\ang{Qz,x-y}\le0$, and
        \item[c.] either $\frac{\sum_{i=1}^n (x_i-y_i)}{\sum_{i=1}^nz_i} \ge 0$, or $\sum_{i=1}^n (x_i-y_i)= {\sum_{i=1}^nz_i} = 0$.
    \end{enumerate}
    Then, $h(x,y)\ge f(x)$.
\end{lemma}
\begin{proof}
    Similar to Lemma~\ref{lemma:direction_concavity}, $f(x)\le h(x,y)$ is equivalent to $\ang{Q(x-y),x-y}\le0$.
    Let $u:=x-y$, and choose $w\in \R^n$ such that 
    \begin{equation*}
        w:=\al z,\quad \text{ where } \al:= 
        \begin{cases} 
            1 &\text{ if } \sum_{i=1}^nz_i = 0,\\
            \frac{\sum_{i=1}^n u_i}{\sum_{i=1}^nz_i}&\text{ otherwise.}
        \end{cases}
    \end{equation*}
    From (c), $\al \ge 0$ and $\sum_{i=1}^nu_i = \sum_{i=1}^nw_i$, or equivalently $\sum_{i=1}^n(u_i-w_i)=0$. Note that from (a) and (b), we have
    \begin{equation*}
        \ang{Qw,w}=\al^2\ang{Qz,z} \ge 0,\quad \ang{Qw,u}\le 0.
    \end{equation*}
    Because $Q$ is conditionally negative definite, we have $\ang{Q(w-u),w-u} \le 0$.
    Combining this with $\ang{Qw,w} \ge 0$ and $\ang{Qw,u}\le 0$, we get
    \begin{equation*}
        \ang{Qu,u} = \ang{Q(w-u),w-u} - \ang{Qw,w} + 2 \ang{Qw,u}\le 0,
    \end{equation*}
    therefore we have that $f(x)\le h(x,y)$.
\end{proof}
Using Lemmas~\ref{lemma:direction_concavity} and~\ref{lemma:setup_valid_tangents}, we now establish conditions for when a tangent plane $h(x,y)$ provides an upper approximation for higher value solutions in $\mathcal{K}$, i.e., $f(x)\ge f(y) \implies h(x,y)\ge f(x)$.
\begin{theorem}
    \label{thm:valid_tangents}
    Suppose $x,y\in\R^n_+$, such that $f(x)\ge f(y)$. Then, $h(x,y)\ge f(x)$ if either
    \begin{enumerate}
        \item[a.] $\sum_{i=1}^n x_i \le \sum_{i=1}^n y_i$, or
        \item[b.] there is $w\in \R^n_{+}\setminus\{0\}$ such that $\ang{Qw,x-y} \le 0$.
    \end{enumerate}
\end{theorem}
\begin{proof}
    Because $f(x) \ge f(y)$, we have
    \begin{equation*}
        \ang{Q(x+y),x-y} \ge 0.
    \end{equation*}
    \begin{enumerate}
        \item[a.] 
        Suppose $\sum_{i=1}^n x_i \le \sum_{i=1}^n y_i$. 
        Choose $z:= -(x+y)$. 
        Because $Q$ has positive off-diagonal entries, then $$\ang{Qz,z} = \ang{-Q(x+y),-(x+y)} =\ang{Q(x+y),x+y} \ge 0,$$ and $\ang{Qz,x-y} \le 0$. 
        Taking into account that $\sum_{i=1}^n x_i \le \sum_{i=1}^n y_i$ and $x+y\in \R^n_+$, condition (c) in Lemma~\ref{lemma:setup_valid_tangents} is fulfilled. 
        Hence, by Lemma~\ref{lemma:setup_valid_tangents}, the inequality $h(x,y)\ge f(x)$ holds.
        \item[b.] 
        Suppose there is $w\in \R^n_{+}\setminus\{0\}$ such that $\ang{Qw,x-y} \le 0$. 
        Then, there is $z \in [w,x^*+y]\subset \R_+^n\setminus\{0\}$ such that $\ang{z,x^*-y} =0$, meaning there is $t\in [0,1]$ such that $$z = t(x^*+y)+(1-t)w \in \R_+^n\setminus\{0\},$$ and $\ang{z,x^*-y} = t\ang{x^*+y,x^*-y}+ (1-t)\ang{w,x^*-y} = 0$. 
        Note that $Q$ has zero diagonal and positive off-diagonal entries, hence $\ang{Qz,z} >0$.
        Therefore by Lemma~\ref{lemma:direction_concavity}, we have that the inequality $h(x,y)\ge f(x)$ holds.
    \end{enumerate}
\end{proof}

\subsection{Cutting plane algorithms}
\label{subsec:algorithms}

We now introduce two cutting plane algorithms designed to solve the nonconcave quadratic problem~\eqref{prob:emsp}.
Let $A\subset \R^n_+$ denote a set of valid tangent planes, where for all $y\in A$ we have that $f(x^*)\le h(x^*,y)$.
Then, we define
\begin{equation*}
    \Gamma_A :=\left\{(x,\theta)\in \R^{n+1}:\; x\in \mathcal{K},\; \theta \le h(x,y),\; \forall y\in A\right\}.
\end{equation*}
The cutting plane model of the~\eqref{prob:emsp} is then given as the following mixed-integer linear program,
\begin{equation}
    \label{prob:ILP_A}
    \tag{$\text{ILP}_{A}$}
    \max_{(x,\theta)\in \Gamma_A}\theta.
\end{equation} 
Given $A$ is a set of valid cuts,~\eqref{prob:ILP_A} provides a valid upper approximation of~\eqref{prob:emsp}.
We now present two algorithms for solving the~\eqref{prob:emsp} that iteratively generate new, valid tangent planes, thereby tightening the approximation of~\eqref{prob:ILP_A}.
Provided the first cut added is valid, both methods are guaranteed to converge to an optimal solution of the~\eqref{prob:emsp}.
Note that from Theorem~\ref{thm:valid_tangents}.a, we can always choose the first cut to be the solution of the maximum cardinality problem, $\max_{x\in \mathcal{K}}\sum_{i=1}^n x_i$.

The first algorithm makes use of the following proposition, which asserts that the tangent plane at the optimal solution of~\eqref{prob:ILP_A} is always valid.
\begin{proposition}
    \label{prop:lp_tangent_validity}
    Given $A\subset \mathbb R^n_+$, let $(x,\theta)$ be an optimal solution of the cutting plane problem~\eqref{prob:ILP_A}.
    Then, there is a $y\in A$ such that $\ang{Qy,x^*-x}\le0$, and hence $f(x^*) \le h(x^*,x)$, where $x^*$ is an optimal solution of~\eqref{prob:emsp}.
\end{proposition}
\begin{proof}
    Suppose, for a contradiction, that for all $y\in A$ we always have $\ang{Qy,x^*-x} >0$, or equivalently, $\ang{Qy,x^*} > \ang{Qy,x}$.
    Then,
    \begin{equation*}
        \theta\le \ang{Qy,x - y} + \tfrac{1}{2}\ang{Qy,y}< \ang{Qy,x^* - y} + \tfrac{1}{2}\ang{Qy,y}
    \end{equation*}
    holds for all $y\in A$.
    Let $\hat\theta$ be such that,
    \begin{equation*}
        \hat \theta := \min_{y \in A} \ang{Qy,x^* - y} + \tfrac{1}{2}\ang{Qy,y}>\theta.
    \end{equation*}
    However, $(x^*,\hat \theta)\in \Gamma_{A}$, and $\hat{\theta} > \theta$, which contradicts $(x,\theta)$ being an optimal solution.
    Hence, the first assertion is settled. 
    The second assertion is a direct consequence of Theorem~\ref{thm:valid_tangents}.b, where $w=y$, and noting that $f(x)\le f(x^*)$. 
    Hence, $f(x^*) \le h(x^*,x)$.
\end{proof}
Using this result, we can now solve the~\eqref{prob:emsp} by repeatedly solving~\eqref{prob:ILP_A} to optimality, and use the solutions as a new valid tangent planes.
An implementation of this approach is shown in Algorithm~\ref{algo:repeated_ilp}, and its convergence is established in Proposition~\ref{prop:algo1_convegence}.

\begin{algorithm}[H]
    \SetAlgoLined
    \vspace{0.5em}
    $k \gets 0$,
    $UB_k\gets +\infty$ \\
    Take $x^0\in \arg\max_{x\in \mathcal{K}} \sum_{i=1}^n x_i$\\
    Set $A_1  \gets \{x^0\}$, $LB_k \gets f(x^k)$ \\
    \While{  $UB_k >LB_k$}{
        $k \gets k+1$ \\        
        Solve $(\text{ILP}_{A_k})$ to obtain  $(x^{k},\theta^k)$ \\
        $UB_k\gets \theta_k$, $LB_k\gets \max\{LB_{k-1},f(x^k)\}$\\
        $A_{k+1} \gets A_{k}\cup\{x^k\}$         
    }
    \label{algo:repeated_ilp}
    \caption{Repeated~\eqref{prob:ILP_A} method for solving~\eqref{prob:emsp}.}
\end{algorithm}

The repeated~\eqref{prob:ILP_A} algorithm is similar to the extended cutting plane method presented in~\cite{westerlund1995extended}, with a modification on the first cut added $h(x,x^0)$. 

\begin{proposition}
    \label{prop:algo1_convegence}
    Algorithm~\ref{algo:repeated_ilp} convergences to an optimal solution of the~\eqref{prob:emsp} in a finite number of steps.
\end{proposition}
\begin{proof}
    As every $(\text{ILP}_{A_k})$ is solved to optimality, we have from Proposition~\ref{prop:lp_tangent_validity} that the tangent plane of every $x^k$ is valid.
    This implies that $(x^*,f(x^*))$ is always feasible at every step $k$, i.e., $(x^*,f(x^*))\in \Gamma_{A_k}$ for all $k\ge 0$. 
    Thus,
    \begin{equation*}
        \text{UB}_k = \max_{(x,\theta)\in \Gamma_{A_k}} \theta \ge f(x^*) =\max_{x\in \mathcal K} f(x) \ge \text{LB}_k.
    \end{equation*}
    Because the feasible region $\mathcal{K}$ is finite (variables $x$ are discrete, and the polyhedral set $P$ is bounded), there is a step $k$ such that the optimal solution $(x^k,\theta^k)$ of $(\text{ILP}_{A_k})$ is such that $x^k \in A_k$. 
    In this case, we have $\text{UB}_k = \theta^k \le h(x^k,x^k) = f(x^k) \le \text{LB}_k$, and hence, $\text{UB}_k =\text{LB}_k$.
    When $\text{UB}_k =  \text{LB}_k$, 	we have $\theta^k = \max_{x\in \mathcal{K}} f(x)$, and therefore Algorithm~\ref{algo:repeated_ilp} converges to an optimal solution. 
\end{proof}

While Algorithm~\ref{algo:repeated_ilp} is globally convergent, it requires solving~\eqref{prob:ILP_A} to optimality at every iteration.
Depending on $\mathcal{K}$, this potentially represents a difficult mixed-integer program.
Overcoming the difficulty of repeatedly solving~\eqref{prob:ILP_A} is usually achieved through a branch and cut implementation, where cuts are added on the fly during the solve procedure.
However, in the case of the~\eqref{prob:emsp}, feasible solutions are not always guaranteed to be valid cuts.
In Algorithm~\ref{algo:forced_cardinality}, we ensure feasible solutions provide valid cuts by iteratively imposing cardinality constraints.
This begins by solving~\eqref{prob:ILP_A} with cardinality forced at its maximum.
This subproblem can then be solved using a branch and cut methodology, and from Theorem~\ref{thm:valid_tangents}.a, any feasible solution forms a valid tangent plane.
We then solve for an upper bound of all future iterations with smaller cardinality.
Then cardinality is then decreased by 1, and the subproblem resolved until the upper bound of lower cardinality problems drops below the best lower bound.
This procedure is outlined in Algorithm~\ref{algo:forced_cardinality}, and convergence is established in  Proposition~\ref{prop:algo2_convegence}.

\begin{algorithm}[H]
    \SetAlgoLined
    \vspace{0.5em}
    $k \gets 0$, $UB_0\gets +\infty$\\
    Take $x^0\in \arg\max_{x\in \mathcal{K}} \sum_{i=1}^n x_i$ \\
    $LB_k \gets f(x^k)$ \\
    $c_1 \gets \sum_{i=1}^n x^0_i$, $A_1  \gets \{x^0\}$ \\
    \While{$UB_k>LB_k$}{
        $k \gets k+1$ \\
        Solve $\max_{(x,\theta)\in \Gamma_{A_k}} \left\lbrace \theta : \sum_{i=1}^n x_i = c_k  \right\rbrace$ for $(x^k,\theta^k)$ using branch and cut \\
        Save all cuts found during the branch and cut procedure and add them to $A_{k+1}$ \\
        Solve $\max_{(x,\theta)\in \Gamma_{A_{k+1}}} \{ \theta: \sum_{i=1}^n x_i \le c_k-1\}$ for $UB_k$ \\
        $LB_k\gets \max\{LB_{k-1},f(x^k)\}$\\
        $c_{k+1} \gets c_k - 1$
    }
    \label{algo:forced_cardinality}
    \caption{Forced cardinality method for solving~\eqref{prob:emsp}.}
\end{algorithm}
\begin{proposition}
    \label{prop:algo2_convegence}
    Algorithm~\ref{algo:forced_cardinality} convergences to the optimal solution of the~\eqref{prob:emsp} in a finite number of steps.
\end{proposition}
\begin{proof}
    From Theorem~\ref{thm:valid_tangents}.a, tangent planes are always valid on solutions that have the same cardinality.
    Hence, at iteration $k$, solving the subproblem $\max_{(x,\theta)\in \Gamma_{A_k}} \left\lbrace \theta : \sum_{i=1}^n x_i = c_k  \right\rbrace$ using a branch and cut methodology gives an optimal solution to the problem $\max_{k\in \mathcal{K}}\left\lbrace f(x) : \sum_{i=1}^n x_i = c_k \right\rbrace $.
    Futhermore, from Theorem~\ref{thm:valid_tangents}.a, all tangent planes remain valid for future iterations with smaller cardinality, and hence
    \begin{equation*}
        UB_k = \max_{(x,\theta)\in \Gamma_{A_{k+1}}} \left\lbrace \theta: \sum_{i=1}^n x_i \le c_k-1\right\rbrace \ge \max_{k\in \mathcal{K}}\left\lbrace f(x) : \sum_{i=1}^n x_i \le c_k-1 \right\rbrace.
    \end{equation*}
    Moveover,
    \begin{equation*}
        LB_k = \max_{k\in \mathcal{K}}\left\lbrace f(x) : \sum_{i=1}^n x_i \ge c_k \right\rbrace 
    \end{equation*}
    and hence if $UB_k \le LB_k$, then $\max_{k\in \mathcal{K}}\left\lbrace f(x) : \sum_{i=1}^n x_i \le c_k-1 \right\rbrace \le \max_{k\in \mathcal{K}}\left\lbrace f(x) : \sum_{i=1}^n x_i \ge c_k \right\rbrace$ and hence the optimal solution has already been found.
    This and the finite domain of the feasible set $\mathcal{K}$ guarantee the convergence of the algorithm.
\end{proof}

In difficult instances of the~\eqref{prob:emsp}, a large number of tangent planes are potentially required to sufficiently approximate the objective function (such as with high-coordinate instances in~\cite{spiers2023}).
Building a large set of strong tangent planes may take many iterations, especially in the case of Algorithm~\ref{algo:repeated_ilp}.
In order to speed up cut generation, observe that Proposition~\ref{prop:lp_tangent_validity} still holds when the integrality of~\eqref{prob:ILP_A} is relaxed.
Therefore, provided $A$ contains valid cuts, the solution of the continuous relaxation of~\eqref{prob:ILP_A} is also guaranteed to provide a valid tangent plane.
The process of generating cuts from the continuous relaxation can be done quickly and is shown in Algorithm~\ref{algo:lp_tangents}.
These cuts can then be added at any stage of the previous two Algorithms.

\begin{algorithm}[H]
    \SetAlgoLined
    \vspace{0.5em}
    $k \gets 0$, $UB_0\gets +\infty$, $LB_0 \gets 0$\\
    Assume $A_0$ is already populated with valid tangent planes.\\
    \While{$UB_k>LB_k$}{
        $k \gets k+1$ \\        
        Solve the continuous relaxation of $(\text{ILP}_{A_k})$ to obtain  $(x^{k},\theta^k)$ \\
        $UB_k\gets \theta_k$, $LB_k\gets \max\{LB_{k-1},f(x^k)\}$, $A_{k+1} \gets A_{k}\cup\{x^k\}$       
    }
    \label{algo:lp_tangents}
    \caption{LP-relaxation cuts for~\eqref{prob:emsp}.}
\end{algorithm}

\section{Numerical results}
\label{sec:results}

We now present numerical results for cutting plane Algorithms~\ref{algo:repeated_ilp} and~\ref{algo:forced_cardinality}. 
These algorithms were implemented in Python 3.10.12 using \texttt{CPLEX} version 22.1.0 as its mixed-integer linear solver. 
The branch and cut method in Algorithm~\ref{algo:forced_cardinality} utilized the \emph{lazy constraint callback} function, enabling the addition of tangent planes as constraints during the branch and bound procedure. 
For each algorithm, LP-relaxation tangents generated by Algorithm~\ref{algo:lp_tangents} are incorporated either at the start of the first iteration ($k=0$), every iteration ($k\ge0$), or not at all. 
This results in six distinct solver configurations. 
Our implementation's source code can be accessed at \url{https://github.com/sandyspiers/euclidean\_maximisation}. 
All tests were conducted on a machine with a 2.3 GHz AMD EPYC processor with 32GB RAM, using a single thread.

The performance of the algorithms was evaluated against the well-known Glover linearisation of the objective function.
This reformulation was first introduced in~\cite{glover_improved_1975}, and is given as
\begin{alignat}{3}
    \label{prob:glover}
    \max \quad & \sum_{i=1}^{n-1} w_{i}, && \\
    \nonumber
    \text{s.t.} \quad
    & x\in P\cap \{0,1\}^n, &&\\
    \nonumber
    & w_{i} \le x_i \sum_{j=i+1}^n q_{ij}, &\quad& 1 \le i \le n-1,\\
    \nonumber 
    & w_{i} \le \sum_{j=i+1}^n q_{ij}x_j, &\quad& 1 \le i \le n-1,\\
    \nonumber 
    & w_{i} \ge 0, &\quad& 1 \le i \le n-1.
\end{alignat}
This formulation was shown in~\cite{Marti2010} to be effective for diversity-sum problems, and was later used as the exact solver for the comprehensive empirical analyses presented in~\cite{parreno2021measuring} and~\cite{Marti2022}.
In addition to~\eqref{prob:glover}, we solve the~\eqref{prob:emsp} using the mixed-integer quadratic programming solver available within \texttt{CPLEX}.

\subsection{Capacitated diversity problem}

We begin by evaluating the performance of the different solution methods for solving the capacitated diversity problem.
In this problem, the constraint set $P$ contains only the following knapsack constraint,
\begin{equation*}
    \sum_{i=1}^n c_i x_i\le b,
\end{equation*}
where $c_i\in \R_+$ ($i=1,\ldots,n$), and $\min_{i=1,\ldots,n}c_i \le b < \sum_{i=1}^n c_i$.
As such, the~\eqref{prob:emsp} then becomes the problem of selecting a subset of predefined locations, each with a weight, to maximize the sum of the pairwise distances, while keeping the total weight less than or equal to a given limit. 
The capacitated diversity problem belongs to the family of diversity problems, which have a wide variety of practical applications, including facility location, social network analysis and ecological conservation \citep{LU2023119856,lai2018solution,peiro2021heuristics}. 

The test instances used are derived from the publically available \texttt{MDPLIB 2.0}\footnote{Available at \url{https://www.uv.es/rmarti/paper/mdp.html}.} test library~\cite{MDPLIB2}.
Within this test library, we use the Euclidean instances of the capacitated diversity problem.
This includes 10 instances each of sizes 50, 150, and 500.
These instances were generated such that the weight of each node was randomly generated in the range $[1,1000]$, with the capacity set to $b = 0.2\sum_{i=1}^nc_i$ and $b = 0.3\sum_{i=1}^nc_i$, making 60 instances in total.

In addition to the previous publically available test sets, we randomly generate some larger instances of the capacitated diversity problem.
These instances are made up of either $1000,1500,2000,2500$ or $3000$ nodes, where each node contains either $2,10$ or $20$ coordinates.
Each coordinate of a location is uniformly randomly generated in the range $[0,100]$.
The weight of each node is uniformly randomly generated in the range $[1,1000]$, and the capacity is set to $b = 0.2\sum_{i=1}^nc_i$ or $b = 0.3\sum_{i=1}^nc_i$.
For every combination of the number of nodes and the number of coordinates, we generate 5 instances, comprising a total of $150$ test instances in total.

\begin{figure}
    \centering
    \includegraphics[width=0.7\linewidth]{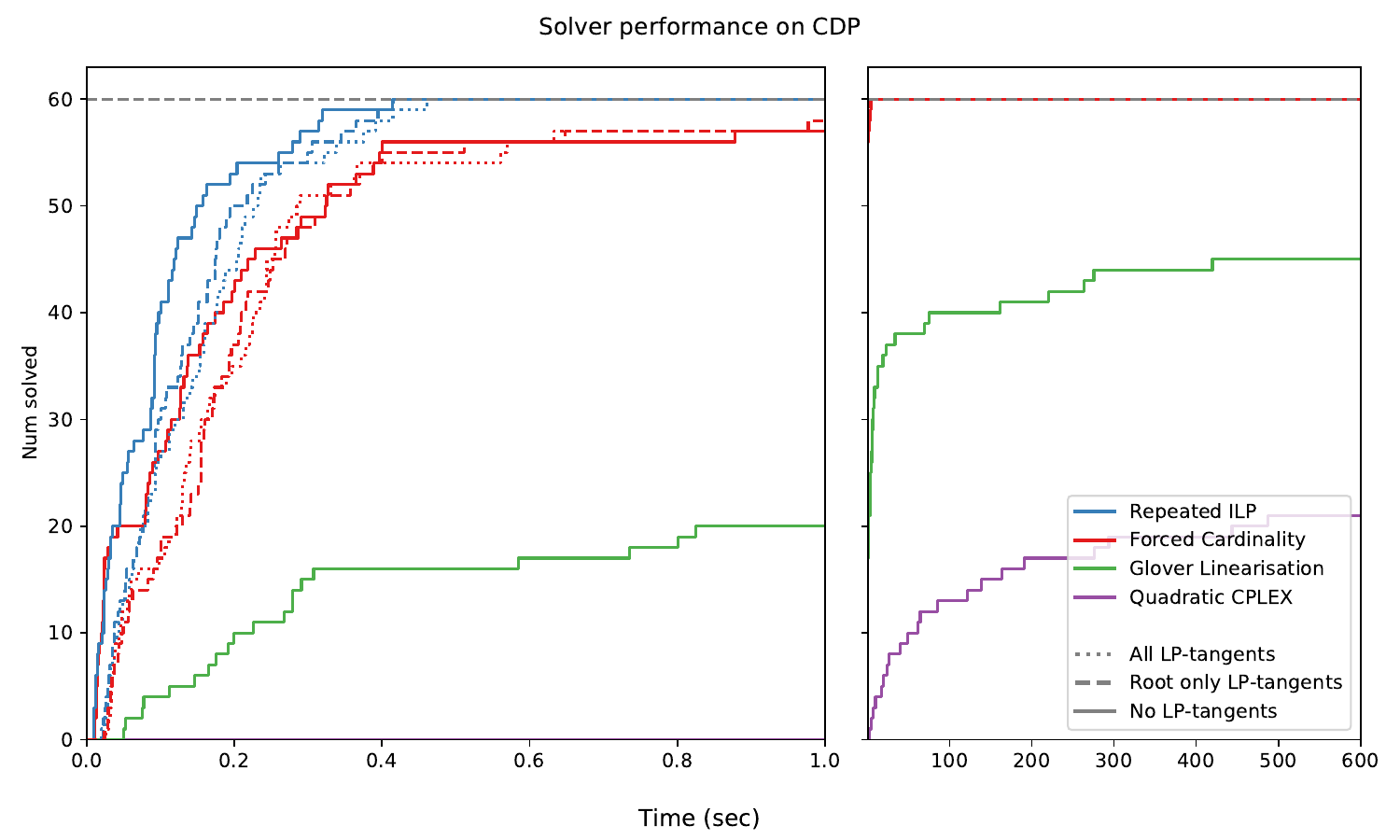}
    \caption{Solver performance on the 60 capacitated diversity problem instances available within the \texttt{MDPLIB 2.0} test library. The repeated~\eqref{prob:ILP_A} and forced cardinality methods are used in three configurations, either adding LP-tangents at the root iteration (dashed line), all iterations (dotted line) or not at all (solid line).  The time axis is split at 1 second due to marked differences in solver performance.}
    \label{fig:pp_cdp}
\end{figure}

\begin{figure}
    \centering
    \includegraphics[width=0.7\linewidth]{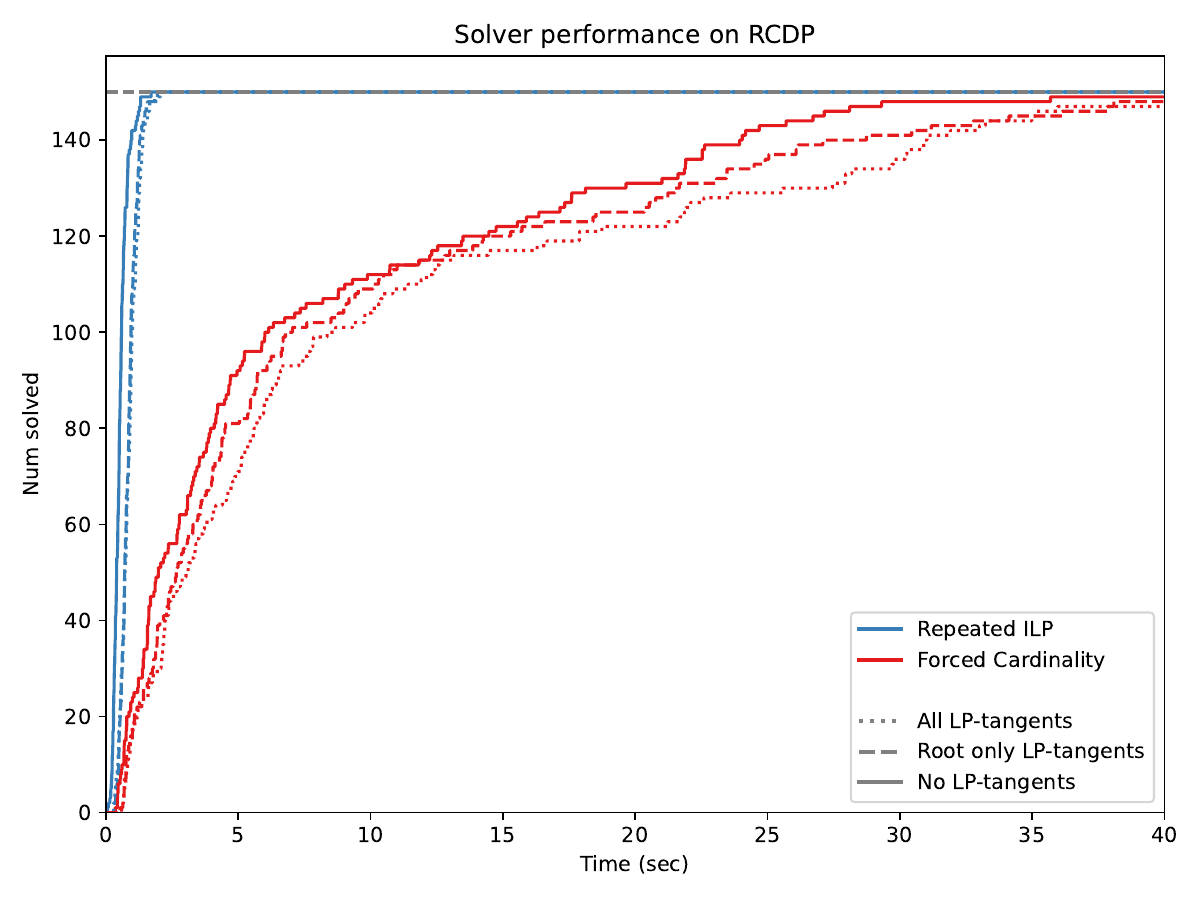}
    \caption{Solver performance on the 150 randomly generated capacitated diversity problem instances. The repeated~\eqref{prob:ILP_A} and forced cardinality methods are used in three configurations, either adding LP-tangents at the root iteration (dashed line), all iterations (dotted line) or not at all (solid line).}
    \label{fig:pp_rcdp}
\end{figure}

The performance of different solver setups for the benchmark problem instances (labelled CDP) and randomized problem instances (labelled RCDP) over a 600-second time limit is displayed in Figures~\ref{fig:pp_cdp} and~\ref{fig:pp_rcdp} respectively. 
As mentioned previously, the repeated~\eqref{prob:ILP_A} and forced cardinality methods are employed in three configurations, either introducing LP-tangents solely at the root iteration, throughout all iterations, or not at all.
Notably, the repeated~\eqref{prob:ILP_A} approach, without LP-tangent planes, emerges as the top performer. 
On CDP test instances, Algorithms~\ref{algo:repeated_ilp} and~\ref{algo:forced_cardinality} exhibit marginal performance differences, both efficiently solving the entire test set within a maximum of 4.25 seconds. 
Furthermore, introducing LP-tangent planes through Algorithm~\ref{algo:lp_tangents} does not enhance performance, and in many cases appears to slow down the process. 
The average solve time for each test set, broken down by problem size, is summarized in Table~\ref{tab:solvetimes}. 
These results highlight significant improvements in the suggested algorithms compared to Glover linearisation and quadratic \texttt{CPLEX}. 
On the RCDP instances, the repeated~\eqref{prob:ILP_A} method substantially outperformed the forced cardinality method.
Remarkably, even with the immense size of these instances, the repeated~\eqref{prob:ILP_A} method was able to solve all instances in under 5 seconds.
Finally, the inclusion of LP-tangent planes appeared to marginally slow solve times.

\begin{table}
    \centering
    \begin{tabular}{@{\extracolsep{4pt}}lrrrrrrrrr@{}}
        \toprule
        &  & \multicolumn{3}{c}{\thead{Repeated~\eqref{prob:ILP_A}}} & \multicolumn{3}{c}{\thead{Forced Cardinality}} &  &  \\
        \cline{3-5}
        \cline{6-8}
        \thead{Type}  & \thead{$n$} & \thead{All} & \thead{Root\\only} & \thead{None} & \thead{All} & \thead{Root\\only} & \thead{None} & \thead{Glover\\Linearisation} & \thead{Quadratic\\\texttt{CPLEX}} \\
        \midrule
        CDP & 50 & 0.0656 & 0.0560 & 0.0349 & 0.0675 & 0.0732 & 0.0256 & 0.3123 & 109.6763 \\
        CDP & 150 & 0.2042 & 0.1851 & 0.1098 & 0.2146 & 0.2206 & 0.1579 & 32.5117 & 600.0009 \\
        CDP & 500 & 0.1622 & 0.1489 & 0.1420 & 0.5430 & 0.5034 & 0.7654 & 500.0613 & 586.6621 \\
        \midrule
        RCDP & 1000 & 0.5851 & 0.5527 & 0.4728 & 3.8039 & 3.4206 & 3.2874 & - & - \\
        RCDP & 1500 & 0.7206 & 0.7472 & 0.4961 & 5.8795 & 5.0636 & 4.9584 & - & - \\
        RCDP & 2000 & 0.9291 & 0.9355 & 0.5296 & 10.5677 & 9.2850 & 8.0693 & - & - \\
        RCDP & 2500 & 1.0937 & 1.0293 & 0.6010 & 14.1290 & 13.6099 & 11.7082 & - & - \\
        RCDP & 3000 & 1.1478 & 1.0350 & 0.5860 & 15.2764 & 12.8296 & 10.2648 & - & - \\
        \midrule
        GDP & 50 & 0.0516 & 0.0493 & 0.0286 & 0.0450 & 0.0440 & 0.0206 & 0.0739 & 0.4484 \\
        GDP & 150 & 0.2306 & 0.2425 & 0.1455 & 0.1760 & 0.1732 & 0.1097 & 1.1257 & 56.3179 \\
        GDP & 500 & 0.2178 & 0.2044 & 0.2100 & 0.6547 & 0.7193 & 0.6191 & 65.8145 & 245.5700 \\
        \midrule
        RGDP & 1000 & 0.5708 & 0.5640 & 0.2388 & 0.7528 & 0.7466 & 0.3704 & - & - \\
        RGDP & 1500 & 0.7483 & 0.8205 & 0.3072 & 0.9375 & 0.9229 & 0.4594 & - & - \\
        RGDP & 2000 & 0.9653 & 0.8686 & 0.3355 & 1.1568 & 1.0983 & 0.5655 & - & - \\
        RGDP & 2500 & 1.0753 & 1.0424 & 0.3872 & 1.3457 & 1.2363 & 0.6157 & - & - \\
        RGDP & 3000 & 1.2021 & 1.1458 & 0.4110 & 1.4407 & 1.2910 & 0.7226 & - & - \\
        \bottomrule
    \end{tabular}
    \caption{Average solve time in seconds of the various solver setups, broken down by test set and test size. Each problem is solved with a time limit of 600 seconds. The repeated~\eqref{prob:ILP_A} and forced cardinality methods are used in three configurations, either adding LP-tangents at all iterations, only the root iteration or not at all.}
    \label{tab:solvetimes}
\end{table}

\subsection{Generalised diversity problem}

The generalised dispersion problem (GDP) represents a fundamental optimisation problem in the fields of facility location, supply chain management, and network design \citep{martinez2021grasp}. 
At its core, the GDP seeks to strategically position a set of facilities on a network to efficiently serve a given demand distribution. 
This entails optimising not only the allocation of facilities to locations but also considering the spread of these facilities. 
The max-sum GDP is given as
\begin{alignat}{3}
	\label{prob:gdp_f}  
    \tag{GDP-f}
	\max \quad & f(x) && \\
	\nonumber
	\text{s.t.} \quad
        \nonumber
	& \sum_{i=1}^n c_i x_i \ge B,\\
        \nonumber
	& \sum_{i=1}^n a_i x_i \le K,\\
        \nonumber
	& x_i\in\{0,1\},\quad i=1,\dots,n,
\end{alignat}
where $c_i$ and $a_i$ represent the capacity and cost of site $i$.
Sites must be chosen such that the minimum demand $B$ is met, and setup cost is kept below the maximum $K$.
The formulation in~\eqref{prob:gdp_f} considers the capacity to be constant if a facility is open.
A more realistic model considers variable setup costs, where extra capacity can be achieved at a given cost, once the facility is open.
The variable cost version of the GDP is given as
\begin{alignat}{3}
	\label{prob:gdp_v}  
    \tag{GDP-v}
	\max \quad & f(x) && \\
	\nonumber
	\text{s.t.} \quad
        \nonumber
	& \sum_{i=1}^n t_i \ge B,\\
        \nonumber
	& \sum_{i=1}^n \left( a_i x_i + b_i t_i \right) \le K,\\
        \nonumber
	& t_i\le c_ix_i,\quad i=1,\dots,n,\\
        \nonumber
	& t_i\in\mathbb{Z}, x_i\in\{0,1\},\quad i=1,\dots,n.
\end{alignat}
We note that~\eqref{prob:gdp_f} and~\eqref{prob:gdp_v} were first introduced in~\cite{martinez2021grasp} where the objective was to maximise the minimum distance, however for our purposes we have changed this objective to maximise the sum of pairwise distances.

For the GDP, we again use the Euclidean test instances available within the \texttt{MDPLIB 2.0} test library on the~\eqref{prob:gdp_v} model.
All parameters were uniformly randomly generated as follows.
The capacity $c_i$ was generated in the range $[1,1000]$, the fixed cost $a_i$ in the range $[c_i/2,2c_i]$ and finally the variable cost $b_i$ in the range $[\min\{1,a_i\}/100,\max\{1,a_i\}/100]$.
The minimum capacity is set at either $B = 0.2\sum_{i=1}^nc_i$ or $B = 0.3\sum_{i=1}^nc_i$.
Finally, the maximum budget is set as $K=\phi \sum_{i=1}^n\left(a_i+b_ic_i\right)$, where $\phi=0.5$ or $\phi=0.6$.
As before, there are 10 instances each of size 50, 150 and 500, making a total of 120 test instances.

To test the solution algorithms at a larger scale, we generate several large instances of~\ref{prob:gdp_v}.
These instances are generated in the same way as described earlier, however, we now increase the number of locations to $1000,1500,2000,2500$ and $3000$ and generate locations with $2,10$ and $20$ sets of coordinates.
For every combination of the number of nodes and the number of coordinates, we generate 5 instances, comprising a total of $300$ test instances in total.

\begin{figure}
    \centering
    \includegraphics[width=0.7\linewidth]{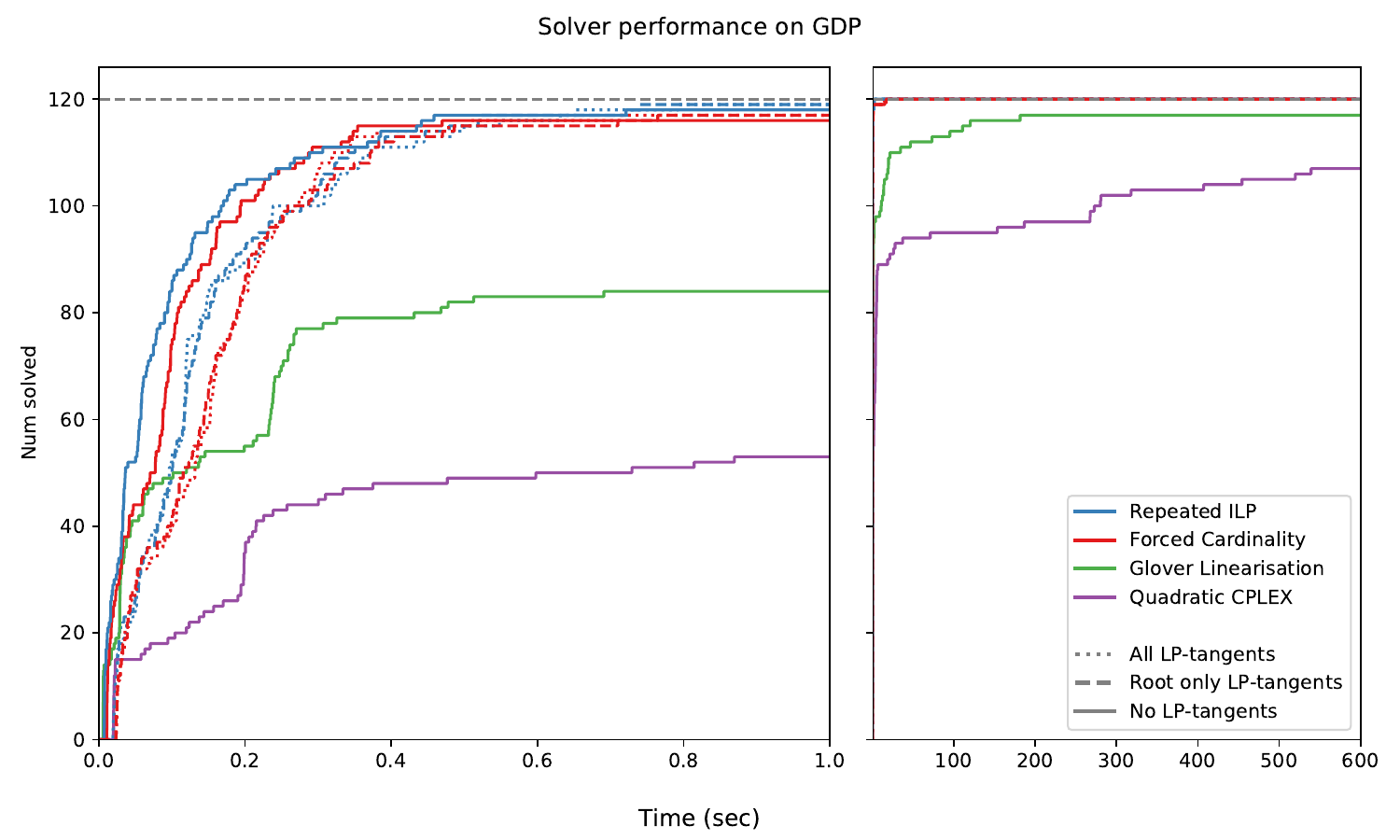}
    \caption{Solver performance on the 120 variable cost generalised diversity problem instances within the \texttt{MDPLIB 2.0} test library. The repeated~\eqref{prob:ILP_A} and forced cardinality methods are used in three configurations, either adding LP-tangents at the root iteration (dashed line), all iterations (dotted line) or not at all (solid line). The time axis is split at 1 second due to marked differences in solver performance.}
    \label{fig:pp_gdp}
\end{figure}

\begin{figure}
    \centering
    \includegraphics[width=0.7\linewidth]{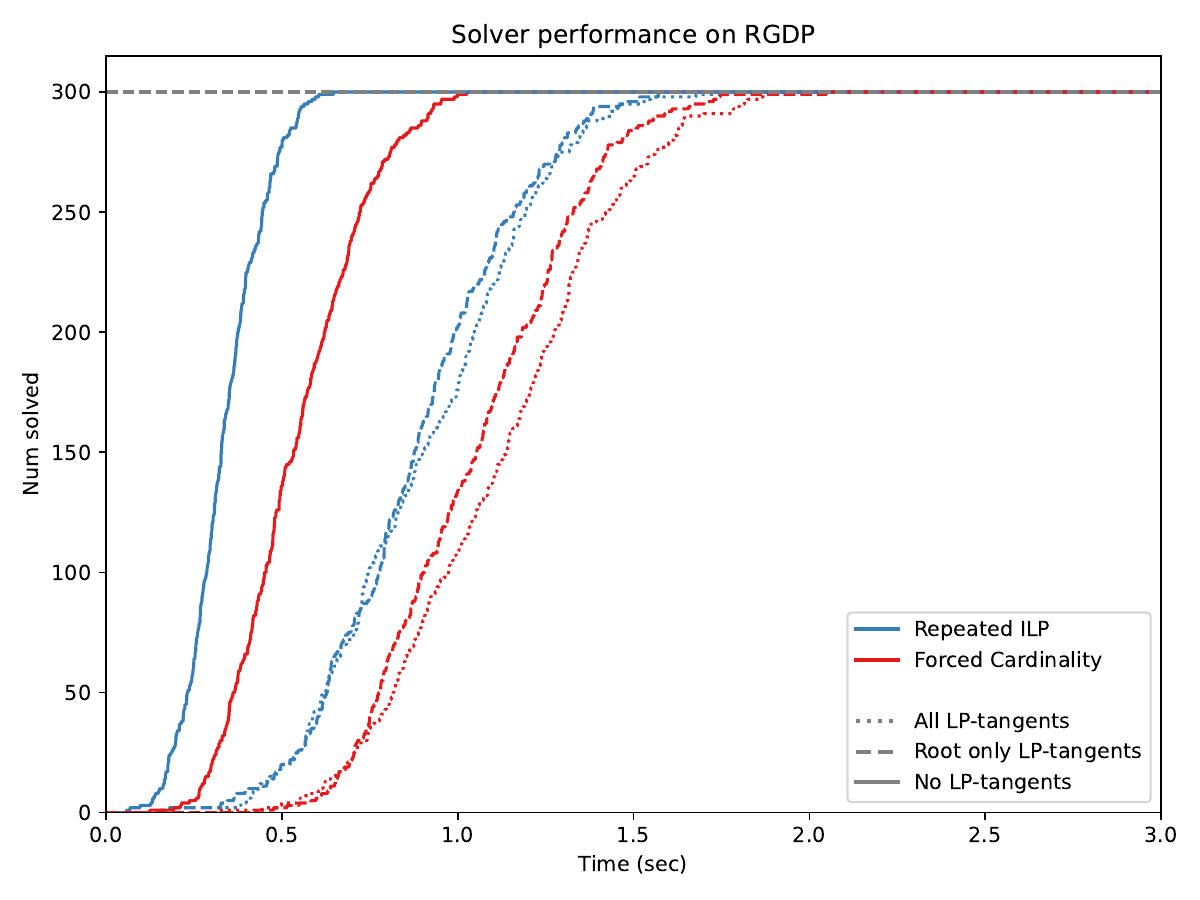}
    \caption{Solver performance on the 300 randomized variable cost generalised diversity problem instances. The repeated~\eqref{prob:ILP_A} and forced cardinality methods are used in three configurations, either adding LP-tangents at the root iteration (dashed line), all iterations (dotted line), or not at all (solid line).}
    \label{fig:pp_rgdp}
\end{figure}

The performance of different solver setups for the benchmark instances (labelled GDP) and random instances (labelled RGDP) over a 600-second time limit is displayed in Figures~\ref{fig:pp_gdp} and~\ref{fig:pp_rgdp} respectively. 
For GDP instances, Algorithms~\ref{algo:repeated_ilp} and~\ref{algo:forced_cardinality} exhibit similar performance, both efficiently solving nearly the entire set within a second. 
The incorporation of LP-tangent planes in either configuration minimally affects the speed of these methods, with differences becoming negligible as solve time surpasses half a second. 
However, both Glover linearisation and Quadratic \texttt{CPLEX} find this test set comparatively easier, as the Glover linearisation model can solve two-thirds of the instances within a second.
Turning to larger randomised GDP instances, as illustrated in Figure~\ref{fig:pp_rgdp}, the repeated~\eqref{prob:ILP_A} method continues to outperform other solver setups. 
Moreover, the results reveal that introducing LP-tangent planes through Algorithm~\ref{algo:lp_tangents} significantly hinders the model's speed. 
A summary of solvetimes for these larger instances is provided in Table~\ref{tab:solvetimes}.

\begin{figure}
    \centering
    \includegraphics[width=0.6\linewidth]{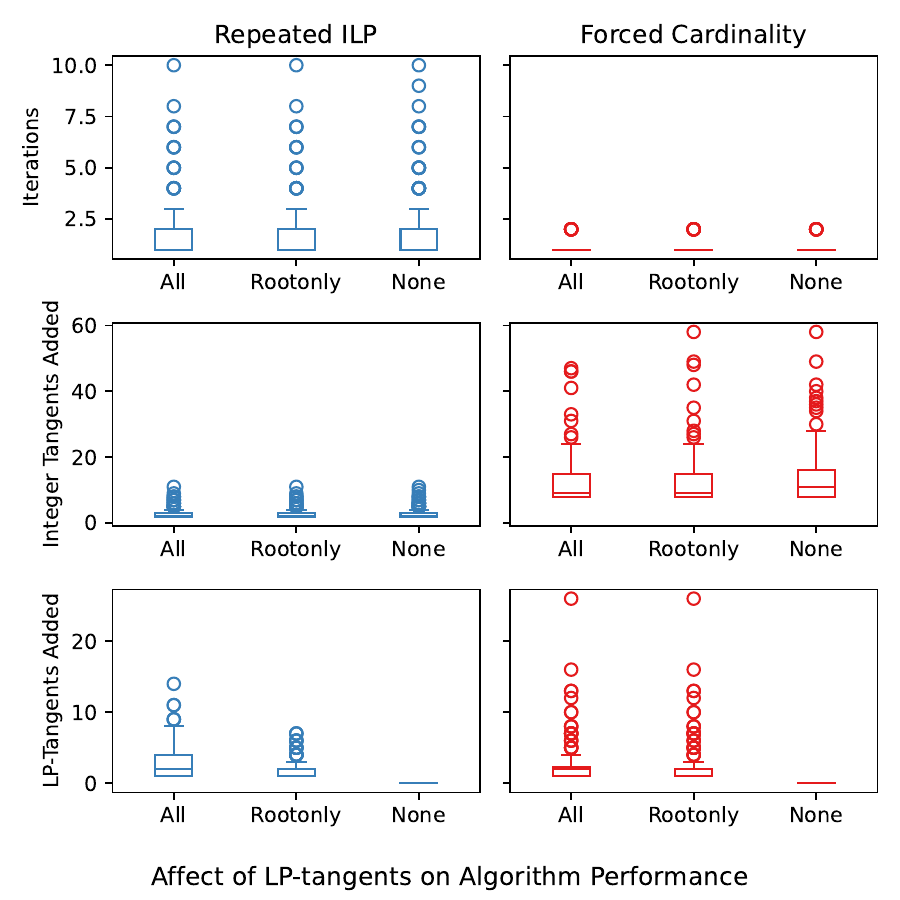}
    \caption{Breakdown of the number of iterations, integer and LP tangents added for each of solver configuration on the CDP and GDP test instances.}
    \label{fig:lp_tan_impact}
\end{figure}

To gain a deeper insight into Algorithms~\ref{algo:repeated_ilp} and~\ref{algo:forced_cardinality}, and to better understand how LP-tangent planes influence their performance, we present a detailed breakdown of each algorithmic setup in Figure~\ref{fig:lp_tan_impact}.
The figure shows the number of iterations and integer- and LP-tangents added across the six solver setups for the CDP and GDP test instances.
Interestingly, the forced cardinality method introduces a significant number of additional integer tangent planes in comparison to the repeated~\eqref{prob:ILP_A} method, despite the latter consistently outperforming in nearly all test sets.
This suggests that by solving~\eqref{prob:ILP_A} to optimality, the cut generated provides a very tight approximation of the objective function at the optimal solution.
Therefore, in many cases, it is worth taking the extra to solve the~\eqref{prob:ILP_A} subproblem to optimality, as the cut generated is expected to be tight.
This also explains why the addition of LP-tangent planes does not seem to provide much computational benefit to either approach.
As these cuts are generated on the continuous relaxation, they are expected to be even further away from the optimal solution than any integer solution, and hence provide a worse approximation.
While LP-tangents are easy to generate and can therefore introduce a large number of cuts, they do not provide a good approximation of the objective function, and hence we see that a similar number of integer tangents are required across all possible LP-tangent configurations.

\subsection{Bi-level diversity problem}

We finish this section with an analysis of the performance of the cutting plane algorithms on the bi-level max-sum diversity problem.
It is noted in~\cite{parreno2021measuring} that the solution of the max-sum diversity problem often contains clusters of nodes close together.
This leads to solutions that may be impractical for many real-world applications.
To overcome this, the authors introduced the bi-level max-sum diversity problem~\eqref{prob:blmsdp}.
The~\eqref{prob:blmsdp} attempts to select a subset of exactly $p$ predefined locations to maximise the sum of pairwise distances, such that all pairwise distances are greater or equal to a given threshold $\delta>0$.
The bi-level diversity problem can be formulated as the following quadratic program,
\begin{alignat}{3}
	\label{prob:blmsdp}  \tag{BLMSDP}
	\max \quad & f(x) && \\
	\nonumber
	\text{s.t.} \quad
        \nonumber
	& \sum_{i=1}^n x_i = p,\\
        \label{prob:blmsdp_threshhold_cts}
        & \sum_{\substack{j=1\\d_{ij} < \delta}}^n x_j \le 1,\quad i=1,\dots,n,\\
        \nonumber
	& x_i\in\{0,1\},\quad i=1,\dots,n.
\end{alignat}
Constraint~\eqref{prob:blmsdp_threshhold_cts} ensures that whenever a node is chosen, none of the nodes that are closer than the threshold are also chosen.
This, therefore, provides a valid formulation of the~\eqref{prob:blmsdp} that fits the structure of the~\eqref{prob:emsp}.

We assessed the performance of the suggested algorithms on the~\eqref{prob:blmsdp} using the location of 3161 Australian postcodes\footnote{The dataset used is available at \url{https://github.com/matthewproctor/australianpostcodes}.}.
Of the 3161 locations, the model was formulated such that $p=50$, and where the Euclidean distance was calculated using the latitude and longitude of a location as its coordinates.
We then solved the model using the six solver setups shown previously, with a time limit of 600 seconds, where the threshold was set to $\delta=0,0.5,1,1.5$ and $2$.

\begin{figure}
    \centering
    \hfill
    \subfigure[$\delta=0$]{
        \includegraphics[width=0.3\textwidth]{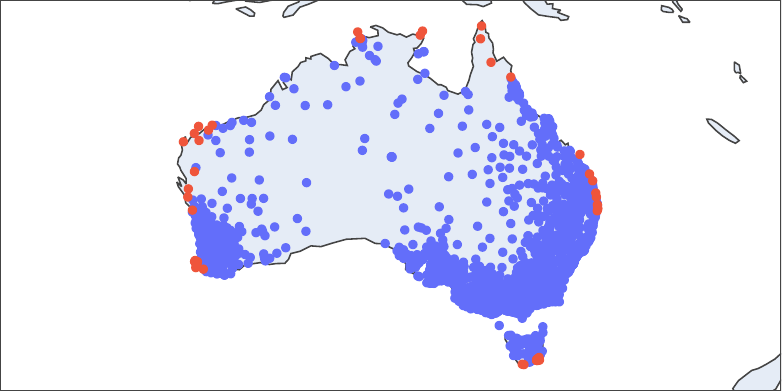}
        \label{fig:figure_a}
    }
    \subfigure[$\delta=0.5$]{
        \includegraphics[width=0.3\textwidth]{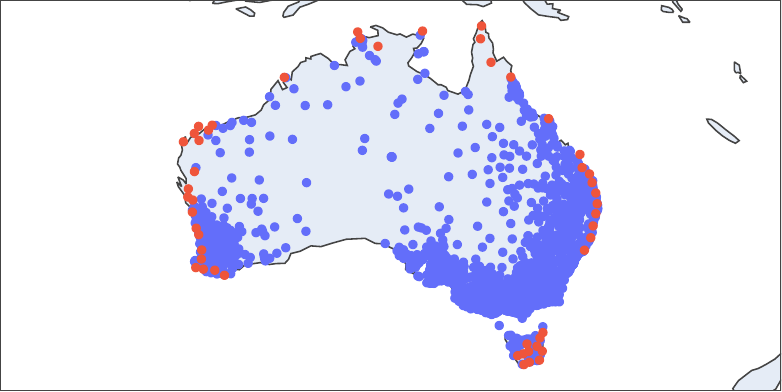}
        \label{fig:figure_b}
    }
    \subfigure[$\delta=1$]{
        \includegraphics[width=0.3\textwidth]{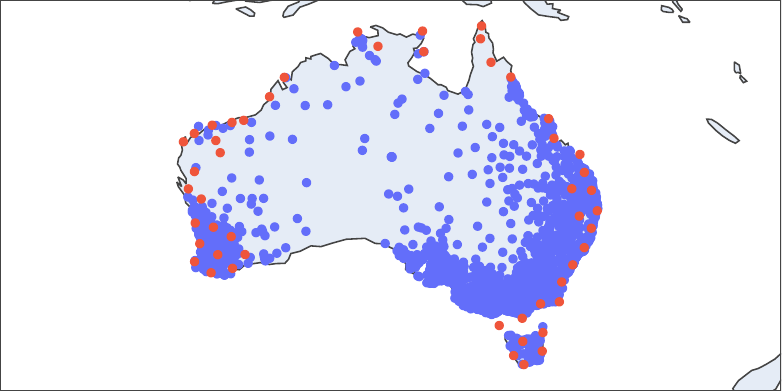}
        \label{fig:figure_c}
    }
    \hfill
    \subfigure[$\delta=1.5$]{
        \includegraphics[width=0.3\textwidth]{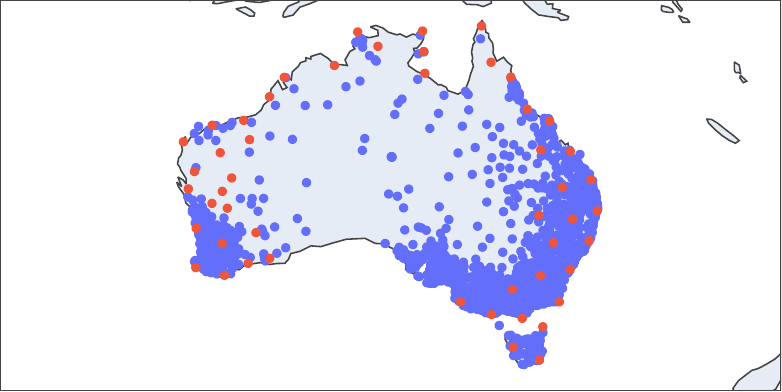}
        \label{fig:figure_d}
    }
    \subfigure[$\delta=2$]{
        \includegraphics[width=0.3\textwidth]{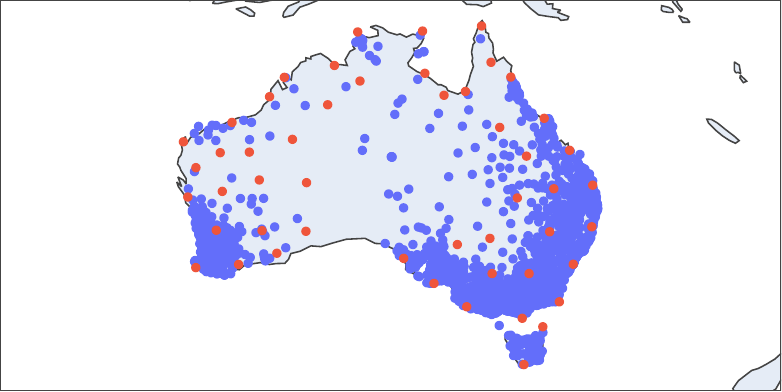}
        \label{fig:figure_e}
    }
    \hfill
    \caption{Visualisation of the optimal solution of the~\eqref{prob:blmsdp} on 3161 Australia postcodes at five minimum distance thresholds.}
    \label{fig:postcodes}
\end{figure}

In Table~\ref{tab:postcodes} we show the solve time and number of cuts required to solve the bi-level problem at each threshold, with the optimal solutions shown in Figure~\ref{fig:postcodes}.
The results appear to contrast those shown in earlier experiments.
It would appear that, in general, the bi-level diversity problem is more difficult than both the CDP and GDP problems.
When $\delta=0$, the problem reduces to the standard Euclidean max-sum diversity problem and appears to be one of the most difficult tested, with both Algorithms~\ref{algo:repeated_ilp} and~\ref{algo:forced_cardinality} reaching the 600-second time-limit when used without LP-tangents.
Furthermore, the addition of LP-tangents appears very useful for these problem instances, especially for low threshold levels.
For instance, when $\delta=0$, adding LP-tangents at every iteration of Algorithm~\ref{algo:repeated_ilp} allowed the model to solve to optimality in under two seconds, whereas without any LP-tangents, the algorithm reached the time-limit without proving optimality.

\begin{table}
    \centering
    \begin{tabular}{@{\extracolsep{4pt}}lrrrrrrr@{}}
        \toprule
        &  &  \multicolumn{3}{c}{\thead{Repeated~\eqref{prob:ILP_A}}} & \multicolumn{3}{c}{\thead{Forced Cardinality}}  \\
         \cline{3-5}
         \cline{6-8}
         & \thead{$\delta$}&\thead{All} & \thead{Root\\only} & \thead{None} & \thead{All} & \thead{Root\\only} & \thead{None} \\
        \midrule
        Solvetime &0&  \textbf{1.9787} & 3.0748 & 600.0159 & 65.4077 & 66.6478 & 600.0402 \\
        Cuts added & &  4 & 6 & 10 & 49 & 49 & 111 \\
         \midrule
        Solvetime &0.5& \textbf{3.6844} & 9.2738 & 33.0451 & 5.3783 & 5.0127 & 17.7160 \\
        Cuts added && 3 & 8 & 26 & 72 & 72 & 128 \\
         \midrule
        Solvetime &1& 9.9957 & 9.5729 & 44.7316 & 8.3350 & \textbf{8.3167} & 13.7821 \\
        Cuts added && 6 & 6 & 20 & 35 & 35 & 79 \\
         \midrule
        Solvetime &1.5& 11.2467 & 10.8825 & 36.7027 & 6.8967 & \textbf{6.8312} & 8.0908 \\
        Cuts added && 5 & 5 & 16 & 30 & 30 & 63 \\
         \midrule
        Solvetime &2& 19.1583 & 22.2657 & 35.2323 & 11.8013 & 11.6925 & \textbf{11.0655} \\
        Cuts added && 6 & 7 & 11 & 31 & 31 & 39 \\
        \bottomrule
    \end{tabular}
    \caption{Solve time in seconds and the number of integer cutting planes added for the six suggested solution methods for the~\eqref{prob:blmsdp} at various threshold levels.  At each threshold, we highlight the fastest solver.}
    \label{tab:postcodes}
\end{table}

\section{Conclusion and future work}

In this paper, we present two exact cutting plane algorithms for the general Euclidean distance maximisation problem.
We establish the validity of tangents by introducing the concept of \emph{directional concavity}.
This notion led to the formulation of two important sufficient conditions for valid cuts, shown in Theorem~\ref{thm:valid_tangents}.
Two cutting plane solution algorithms were then introduced.
The algorithms exploit Theorem~\ref{thm:valid_tangents} to ensure the search for the optimal solution always stays on a concave direction of the objective function, therefore ensuring all cuts are valid.
This was achieved by either repeatedly solving the cutting problem subproblem to optimality, or by iteratively forcing and decreasing the cardinality of the problem.
Furthermore, we showed how cuts can be quickly generated by solving the continuous relaxation of the cutting plane subproblem.

Extensive numerical experiments were used to test the suggested solution algorithms.
The results are very promising, with all proposed methods easily able to solve capacitated diversity problem instances with 3000 locations in under 60 seconds.
This represents a significant improvement compared to other exact methods for the~\eqref{prob:emsp}.
Furthermore, we use the bi-level diversity problem to show how the structure of a problem can change the relative performance of the different approaches.
For instance, for large-scale generalised diversity problems with variable costs, the repeated~\eqref{prob:ILP_A} method without LP-tangents was by far the best performer.
However, for all thresholds of the bi-level problem, this setup performed the worst.
Therefore, the choice of which approach to use should depend on the specific problem structure.

The identification of specific problem structures remains an important avenue for future research.
We note a significant gap in the literature on the application of the~\eqref{prob:emsp} to real-world problems.
While the tests used here provide interesting conceptual frameworks, they have scarcely been applied to real-world datasets and problems.
In addition to identifying practical~\eqref{prob:emsp} models, we should also attempt to identify difficult instances of these problems.
In~\cite{spiers2023} we showed how the diversity problem becomes more challenging with a larger number of coordinates.
That difficulty was not observed for the CDP or GDP problems, and hence more work is required to identify other difficult instances of the~\eqref{prob:emsp}.
These problems can also help to understand and decide on which algorithm to use in which scenario.

\section*{Acknowledgment}
The authors are supported by the Australian Research Council through the Centre for Transforming Maintenance through Data Science (grant number IC180100030). 
This work was supported by resources provided by the Pawsey Supercomputing Research Centre with funding from the Australian Government and the Government of Western Australia.

\bibliography{biblio}

\end{document}